\documentclass[12pt]{article}
\usepackage[T1]{fontenc}
\usepackage[utf8]{inputenc}

\usepackage{siunitx}

\usepackage{amsthm,amsmath,amssymb}

\usepackage[]{hyperref}

\newtheorem{teorema}{Theorem}
\newtheorem{teo}{Theorem}[section]

\theoremstyle{definition}
\newtheorem{ese}[teo]{Example}

\numberwithin{equation}{section}

\newcommand{\R}{\mathbb{R}}

\usepackage{xcolor} 

\usepackage{tikz,xcolor,hyperref}
\definecolor{lime}{HTML}{A6CE39}
\DeclareRobustCommand{\orcidicon}{%
	\begin{tikzpicture}
	\draw[lime, fill=lime] (0,0) 
	circle [radius=0.16] 
	node[white] {{\fontfamily{qag}\selectfont \tiny ID}};
	\draw[white, fill=white] (-0.0625,0.095) 
	circle [radius=0.007];
	\end{tikzpicture}
	\hspace{-2mm}
}

\newcommand{\orcidref}[1]{\href{https://orcid.org/#1}{\orcidicon}}

\title{A new class of separable Lagrangian systems generalizing Sawada-Kotera\footnote{This version of the paper corrects some very small mistakes but is essentially the same as the original first version. The final version, stimulated by the Referees, has and introduction with a larger discussion, literature and an extra non polynomial example appeared in \url{https://www.mdpi.com/2673-8716/4/3/26}}}

\author{
Gianluca Gorni \orcidref{0000-0003-1055-4376}\\
Universit\`a di Udine,
Dipartimento di\\ Scienze Matematiche, Informatiche e Fisiche\\
via delle Scienze~208, 33100 Udine, Italy\\
\tt{gianluca.gorni@uniud.it}
\and 
Mattia Scomparin \orcidref{0000-0002-2795-5929}\\ 
Mogliano Veneto 31021, Italy\\
\tt{mattia.scompa@gmail.com}
\and
Gaetano  Zampieri \orcidref{0000-0003-1315-137X}\\
Universit\`a di Verona,
Dipartimento di Informatica\\
strada Le Grazie 15, 37134 Verona, Italy\\
\tt{gaetano.zampieri@univr.it}}

\date{}
\begin{document}
\maketitle

\begin{abstract}
Some characteristics of the Sawada-Kotera Lagrangian system lend themselves to generalization,
 producing a large class of separable Lagrangian systems with two degrees of freedom.
\end{abstract}

{\textbf{Keywords:} Integrable systems\\ Separation of coordinates\\ Generalized H\'enon-Heiles systems}
\section{Main result and related literature}\label{introduction}

\begin{teorema}\label{mainTheorem}
The two smooth Lagrangians in two degrees of freedom
\begin{equation}
  L=\frac{1}{2}(\dot q_1^2+\dot q_2^2)-U(q_1,q_2),\qquad
  \tilde L=\dot q_1\dot q_2-\tilde U(q_1,q_2)
\end{equation}
have the same Lagrange equations and, hence, share the two energy first integrals
\begin{equation}
  E=\frac{1}{2}(\dot q_1^2+\dot q_2^2)+U(q_1,q_2),\qquad
  \tilde E=\dot q_1\dot q_2+\tilde U(q_1,q_2),
\end{equation}
if and only if there exist two smooth functions $f,g$ of one variable such that
\begin{equation}
  U(q_1,q_2)=f(q_1+q_2)+g(q_1-q_2),\qquad
  \tilde U(q_1,q_2)=f(q_1+q_2)-g(q_1-q_2).
\end{equation}
If this happens, the change of variables $(x,y)=(q_1+q_2,q_1-q_2)$ separates the Lagrange equations into
\begin{equation}
  \ddot x=-2f'(x),\qquad
  \ddot y=-2g'(y).
\end{equation}
\end{teorema}

Our motivation was the study of the Lagrangian
\begin{equation}
  L_b=\frac{1}{2}
  (\dot q_1^2+\dot q_2^2)-
  \biggl(\frac{1}{2}(q_1^2+q_2^2)+
  q_1^2q_2+\frac{1}{3}q_2^3\biggr),
\end{equation}
which we call the Sawada-Kotera system. It is integrable, too, by means of the energy and the supplementary first integral 
\begin{equation}\label{SKfi3}
  K=\dot q_1\dot q_2+q_1q_2+\frac{1}{3}q_1^3+q_1q_2^2.
\end{equation} 
This case is also separable in the coordinates $(q_1+q_2,q_1-q_2)$ and its solution can be expressed through elliptic functions. This result was obtained by Aizawa and Saito~\cite{Stefan1972}, but their names did not stick, and  the system is found associated with Sawada and Kotera, and is connected with soliton theory \cite{Ballesteros2010, Fordy1991, Sawada1974}.

The Sawada-Kotera system is special case of the following family of Lagrangian functions, depending on the parameter $b\in\R$:
\begin{equation}
  L_b=\frac{1}{2}
  (\dot q_1^2+\dot q_2^2)-
  \biggl(\frac{1}{2}(q_1^2+q_2^2)+
  q_1^2q_2-\frac{b}{3}q_2^3\biggr).
\end{equation}
See for instance the larger family of Hamiltonian functions (1) 
in \cite{Stefan1984} (not  the one in the Abstract) which we restricted choosing $\omega_1=\omega_2=a=1$.   Since these Lagrangians are autonomous, the associated Lagrange equations have the first integral of \emph{energy}
\begin{equation}
 E_b=\frac{1}{2}(\dot q_1^2+\dot q_2^2)+
 \biggl(\frac{1}{2}(q_1^2+q_2^2)+
 q_1^2q_2-\frac{b}{3}q_2^3\biggr).
\end{equation}

The special case of $b=1$ is the influential H\'enon-Heiles system that was introduced in~\cite{Henon1964} in order to model a Newtonian axially-symmetric galactic system. Such model has no analytic first integral independent of energy, as proved by Ito \cite{Ito}, and its chaotic dynamic behaviour has been extensively studied \cite{Churchill,Tabor1989,Boccaletti2004}. (Note that Wojciechowski calls after H\'enon-Heiles the whole family in his formula (1) quoted above).

By contrast, two other members are known to be integrable. One is the case $b=-6$ studied  by Dorizzi, Grammaticos and Ramani~\cite{Dorizzi1983,Grammaticos1983}, and independently integrated by Wojciechowski~\cite{Stefan1984}.

The case $b=-1$ is precisely the Sawada-Kotera system, which is integrable, too, as we mentioned above. In the 
same spirit of solitons the case $b=-6$ is called Korteweg-De Vries.

Let us mention a recent paper by Sottocornola \cite{Sottocornola2019}, which deals with separation of variables for seven integrable systems related to H\'enon-Heiles, and presents some open questions.

\section{Two Lagrangians for the Sawada-Kotera equations}\label{section:SK}

Consider the \emph{Sawada-Kotera} cubic Lagrangian 
\begin{equation}\label{SK}
  L=\frac{1}{2}(\dot q_1^2+\dot q_2^2)-
  \biggl(\frac{1}{2}(q_1^2+q_2^2)+
  q_1^2q_2+\frac{1}{3}q_2^3\biggr).
\end{equation}
The associated Sawada-Kotera Lagrange equations are 
\begin{equation}\label{SKeqs}
  \ddot q_1=-q_1(1+2 q_2),\qquad 
  \ddot q_2=-q_2-q_1^2-q_2^2.
\end{equation}
The energy first integral is
\begin{equation}\label{eq:E}
  E=\frac{1}{2}(\dot q_1^2+\dot q_2^2)+
  \biggl(\frac{1}{2}(q_1^2+q_2^2)+
  q_1^2q_2+\frac{1}{3}q_2^3\biggr).
\end{equation}
and there is another first integral quadratic in the velocities
\begin{equation}\label{SKfi3bis}
  K =\dot q_1\dot q_2+\Bigl(q_1q_2+\frac{1}{3}q_1^3+
  q_1q_2^2\Bigr).
\end{equation} 

Let us consider also the Lagrangian function
\begin{equation}
  {\tilde L}=\dot q_1\dot q_2-
  \Bigl(q_1q_2+\frac{1}{3}q_1^3+q_1q_2^2\Bigr).
\end{equation}
It happens that
\begin{gather}
  \partial_{q_1}L-\frac{d}{dt}\partial_{\dot q_1}L=
  \partial_{q_2}\tilde L-
  \frac{d}{dt}\partial_{\dot q_2}\tilde L,\\
  \partial_{q_2}L-\frac{d}{dt}\partial_{\dot q_2}L=
  \partial_{q_1}\tilde L-
  \frac{d}{dt}\partial_{\dot q_1}\tilde L,
\end{gather}
so that the Euler-Lagrange equation for the two system are exactly the same:
\begin{equation}
  \ddot q_2=-q_2-q_1^2-q_2^2,\qquad 
  \ddot q_1=-q_1(1+2 q_2).
\end{equation}

The first integral of energy for $\tilde L$ happens to coincide with \eqref{SKfi3bis} above
\begin{equation}
  \tilde E=\dot q\cdot
  \frac{\partial \tilde L}{\partial \dot q}-\tilde L=
  \dot q_1\dot q_2+q_1q_2+\frac{1}{3}q_1^3+q_1q_2^2=
  K
\end{equation} 
(the central dot is the scalar product).
 
\section{Proof of the main theorem}\label{section:GSK}

Suppose we have a Lagrangian of the form
\begin{equation}\label{eq:LL}
  L=\frac{1}{2}(\dot q_1^2+\dot q_2^2)-U(q_1,q_2),
\end{equation}
where the potential $U$ is smooth but it isn't necessarily a polynomial function. The associated Euler-Lagrange equations are
\begin{equation}\label{eqELell}
  \ddot q_1=-\frac{\partial U}{\partial q_1}(q_1,q_2),\qquad
  \ddot q_2=-\frac{\partial U}{\partial q_2}(q_1,q_2).
\end{equation}

Inspired by the Sawada-Kotera results of the previous Section, let us consider a second Lagrangian of the form
\begin{equation}\label{eq:LLTILDE}
  \tilde L=\dot q_1\dot q_2-\tilde U(q_1,q_2),
\end{equation}
whose equations of motions are
\begin{equation}\label{eqELip}
  \ddot q_2=-\frac{\partial \tilde U}{\partial q_1}(q_1,q_2),\qquad
  \ddot q_1=-\frac{\partial \tilde U}{\partial q_2}(q_1,q_2).
\end{equation}
Let us impose that the equations \eqref{eqELell} coincide with those in \eqref{eqELip}:
\begin{equation}\label{conditionsforcoincidence}
  \frac{\partial U}{\partial q_1}=
  \frac{\partial \tilde U}{\partial q_2},\qquad  
  \frac{\partial U}{\partial q_2}=
  \frac{\partial \tilde U}{\partial q_1}.
\end{equation}
We claim that this occurs if and only if 
\begin{align}\label{solutionOfConditionsForCoincidence}
  U(q_1,q_2)={}&
  f(q_1+q_2)+g(q_1-q_2),\\
  \tilde U(q_1,q_2)={}&
  f(q_1+q_2)-g(q_1-q_2),
\end{align}
for some smooth functions~$f,g$. 

To prove this, let us make a change of dependent variables through a \ang{45} rotation in the $(q_1,q_2)$ plane, introducing new variables $r_1,r_2,V,\tilde V$:
\begin{align}
  V(r_1,r_2):={}&U\biggl(\frac{r_1+r_2}{\sqrt{2}},
  \frac{r_2-r_1}{\sqrt{2}}\biggr),
  \\
  \tilde V(r_1,r_2):={}&
  \tilde U\biggl(\frac{r_1+r_2}{\sqrt{2}},
  \frac{r_2-r_1}{\sqrt{2}}\biggr).
\end{align}
that is to say,
\begin{align}
  U(q_1,q_2)={}&V\biggl(\frac{q_1-q_2}{\sqrt{2}},
  \frac{q_1+q_2}{\sqrt{2}}\biggr),\\
  \tilde U(q_1,q_2)={}&\tilde V\biggl(\frac{q_1-q_2}{\sqrt{2}},
  \frac{q_1+q_2}{\sqrt{2}}\biggr),
\end{align}
In terms of the new variables the differential equations~\eqref{conditionsforcoincidence} become
\begin{equation}
  \frac{\partial(\tilde V+V)}{\partial r_2}=0,\qquad
  \frac{\partial(\tilde V-V)}{\partial r_1}=0,
\end{equation}
which mean that there exist one-variable functions $\varphi,\psi$ such that
\begin{equation}
  \tilde V+V=2\varphi(r_1),\qquad
  \tilde V-V=-2\psi(r_2),
\end{equation}
that is,
\begin{equation}
  \tilde V(r_1,r_2)=\varphi(r_1)-\psi(r_2),\qquad
  V(r_1,r_2)=\varphi(r_1)+\psi(r_2),
\end{equation}
which are equivalent to~\eqref{solutionOfConditionsForCoincidence}.

The Euler-Lagrange equations \eqref{eqELell} have the first integral of energy
\begin{equation}
  E=\dot q_1\partial_{\dot q_1}L+
  \dot q_2\partial_{\dot q_2}L-L=
  \frac{1}{2}(\dot q_1^2+\dot q_2^2)+U(q_1,q_2),
\end{equation}
while the first integral of energy for \eqref{eqELip} is
\begin{equation}\label{tildeE}
  \tilde E=
  \dot q_1\partial_{\dot q_1}\tilde L+
  \dot q_2\partial_{\dot q_2}\tilde L-\tilde L=
  \dot q_1\dot q_2+\tilde U(q_1,q_2).
\end{equation}

Whenever conditions \eqref{conditionsforcoincidence} hold, the Euler-Lagrange equations coincide and $E$ and $\tilde E$ are first integrals for both. In this case, in terms of the new variables $(x, y)=(q_1+q_2,q_1-q_2)$ (just a little simpler than $r_1,r_2$ above) the Lagrangian function $L$ becomes
\begin{equation}
  L=
  \frac{1}{4}\dot x^2-f(x)+
  \frac{1}{4}\dot y^2-g(y),
\end{equation}
for which the Euler-Lagrange equations are separated:
\begin{equation}
  \ddot x=-2f'(x),\qquad
  \ddot y=-2g'(y).
\end{equation}
Similarly, the Lagrangian $\tilde L$ becomes
\begin{equation}
  \tilde  L=
  \frac{1}{4}\dot x^2-f(x)-
  \biggl(\frac{1}{4}\dot y^2-g(y)\biggr).
\end{equation}

Given a Lagrangian function $L(t,q,\dot q)$, $q,\dot q\in\R^n$, it is easy to build other Lagrangians which have the same Lagrange equations: simply take an arbitrary smooth scalar function $G(t,q)$, and define
\begin{equation}
  \mathcal{L}(t,q,\dot q):=
  L(t,q,\dot q)+\frac{\partial G}{\partial t}(t,q)+
  \frac{\partial G}{\partial q} (t,q) \cdot \dot q\,.
\end{equation}
We say $L$ and $\mathcal{L}$ are related by a \emph{gauge transform}.  Our Lagrangian functions $L$ in \eqref{eq:LL}, and  $\tilde L$ in \eqref{eq:LLTILDE}, are \emph{not} related by a gauge transform because their difference is quadratic in the velocities.

\section{Examples}
\label{section:App}

Here are some some sample systems that are covered by our results.

\begin{ese}\textbf{(Recovering Sawada-Kotera)}
The functions
\begin{equation}
  f(x)=\frac{1}{4} x^2+\frac{1}{6}x^3,\qquad
  g(y)=\frac{1}{4} y^2-\frac{1}{6}y^3,
\end{equation}
give the Sawada-Kotera $U, \tilde U$:
\begin{gather}
  U(q_1,q_2)=
  \frac{1}{2}(q_1^2+q_2^2)+q_1^2q_2+
  \frac{1}{3}q_2^3=f(q_1+q_2)+g(q_1-q_2),\\
  \tilde U(q_1,q_2)=
  q_1q_2+\frac{1}{3}q_1^3+q_1q_2^2
  =f(q_1+q_2)-g(q_1-q_2).
\end{gather}
Hence, from \eqref{eqELell}, the separated Euler-Lagrange equations for $x=q_1+q_2$, $y=q_1-q_2$ are
\begin{equation}\label{eqELellSK}
  \ddot x=-x-x^2,\qquad
  \ddot y=-y+y^2.
\end{equation}
\end{ese}

\begin{ese}
The class of system to which our result applies is not limited to polynomial potential. Let us consider the following Lagrangian, that is a variation on Calogero's potential  for $n$-bodies constrained on a line with inversely quadratic pair potentials:
\begin{equation}
\label{eq:cal}
  L(q_1,q_2)=\frac{1}{2}(\dot q_1^2+\dot q_2^2)-
  \biggl(\frac{\alpha_f}{(q_2+q_1)^2}+
  \frac{\alpha_g}{(q_2-q_1)^2}\biggr)\,.
\end{equation}
The equations of motion are
\begin{equation}
  \ddot q_1=\frac{2\alpha_f}{(q_2+q_1)^3}-
  \frac{2\alpha_g}{(q_2-q_1)^3},\qquad
  \ddot q_2=\frac{2\alpha_f}{(q_2+q_1)^3}+
  \frac{2\alpha_g}{(q_2-q_1)^3}.
\end{equation}
The functions
\begin{equation}
  f(x)=\frac{\alpha_f}{x^2},\qquad
  g(y)=\frac{\alpha_g}{y^2}\,,
\end{equation}
give the potential in \eqref{eq:cal} and satisfy the hypotheses of Theorem~\ref{mainTheorem}. Hence, from \eqref{eqELell}, the separated Euler-Lagrange equations are
\begin{equation}
  \ddot x = \frac{4\alpha_f}{x^3},\qquad
  \ddot y = \frac{4\alpha_g}{y^3},
\end{equation}
that have simple solutions.

Incidentally, to this system we can apply the following result \cite[Sec.~3]{GZ}: for general Lagrangian functions $\frac{1}{2}m \lVert\dot q\rVert^2-U(q)$, $q\in\R^n$, with $U$ homogeneous of degree~$-2$, we have that for each motion $\lVert q(t) \rVert^2= \frac{2}{m}(E t^2+I t+J)$ where $E$ is the energy and $I,J$ are constants depending on the motion. In the present case $m=1$ and the relation becomes
\begin{gather}
  q_1(t)^2+q_2(t)^2=2(Et^2+I t+J),\\
  x(t)^2=2(E_1t^2+I_1 t+J_1),\qquad  y(t)^2=2(E_2t^2+I_2 t+J_2).
\end{gather}
\end{ese}


\nocite{*}
\providecommand{\bysame}{\leavevmode\hbox to3em{\hrulefill}\thinspace}


\begin{thebibliography}{10}

\bibitem{Stefan1972}
\textsc{Y. Aizawa,  N. Saito}, 
\emph{On the Stability of Isolating Integrals. Effect of the Perturbation in the Potential Function}, 
J. Phys. Soc. Jpn. \textbf{32} (1972), 1636--1640.

\bibitem{Ballesteros2010}
\textsc{A. Ballesteros, A. Blasco}
\emph{Integrable H\'enon-Heiles Hamiltonians: a Poisson algebra approach},
Annals of Physics \textbf{325} (2010), 2787--2799.

\bibitem{Boccaletti2004}
\textsc{D. Boccaletti, G. Pucacco}, 
\emph{Theory of Orbits}, 
Springer, Berlin, 2004.

\bibitem{Churchill}
\textsc{R.C. Churchill, D.L. Rod},
\emph{Geometrical aspects of Ziglin's non-integrability theorem for complex Hamiltonian systems},
J. Diff. Eqns. \textbf{76} (1988), 91--114.

\bibitem{Dorizzi1983}
\textsc{B. Dorizzi, B. Grammaticos, and A. Ramani}, 
\emph{Integrability of Hamiltonians with third and fourth degree polynomial potentials}, 
J. Math. Phys. \textbf{24} (1983), 2289--2295.

\bibitem{Fordy1991}
\textsc{A.P. Fordy}, 
\emph{The H\'enon-Heile system revisited}, 
Phys. D \textbf{52} (1991), 204--210.

\bibitem{GZ}
\textsc{G. Gorni G., G. Zampieri},
\emph{Lagrangian dynamics by nonlocal constants of motion},
Discrete Contin. Dyn. Syst. Ser. S \textbf{13} (2020), no. 10, 2751--2759.

\bibitem{Grammaticos1983}
\textsc{B. Grammaticos, B. Dorizzi, and A. Ramani}, 
\emph{Integrability of Hamiltonians with third and fourth degree polynomial potentials}, 
J. Math. Phys. \textbf{24} (1983), 2289--2295.

\bibitem{Henon1964}
\textsc{M. Henon, C. Heiles}, 
\emph{The applicability of the third integral of motion: Some numerical experiments}, 
Astron. J. \textbf{69} (1964), 73.

\bibitem{Ito}
\textsc{H. Ito}, \emph{Non-integrability of H\'enon-Heiles system and theorem of Ziglin}, Kodai Math. J. \textbf{8} (1985), 120--138.

\bibitem{Ramani1982}
\textsc{A. Ramani, B. Dorizzi and B. Grammaticos}, 
\emph{Painlev\'{e} conjecture revisited}, 
Phys. Rev. Lett. \textbf{24} (1982), 1539--1541.

\bibitem{Sawada1974}
\textsc{K. Sawada and T. Kotera},
\emph{A Method for Finding N-Soliton Solutions of the K.d.V. Equation and K.d.V.-Like Equation},
Prog. Theor. Phys. \textbf{51} (1974), 1355.

\bibitem{Sottocornola2019}
\textsc{N. Sottocornola}, 
\emph{Separation coordinates in H\'enon-Heiles systems}, 
Physics Letters A 383 \textbf{36} (2019), 126027.

\bibitem{Tabor1989}
\textsc{M. Tabor}, 
\emph{Chaos and Integrability in Nonlinear Dynamics}, 
Wiley, New York, 1989.

\bibitem{Stefan1984}
\textsc{S. Wojciechowski}, 
\emph{Separability of an integrable case of the H\'enon-Heiles system}, 
Physics Letters A \textbf{100}, 6 (1984), 277--278.


\end{thebibliography}
\end{document}